\documentclass{article}
\usepackage{amssymb,amsmath,amsthm,longtable}
\usepackage[dvips]{graphicx}
\newcommand{\D}{{\mathbb D}}
\newcommand{\C}{{\mathbb C}}
\newcommand{\E}{{\mathbb E}}

\renewcommand{\Re}{\mathrm {Re\,}}
\renewcommand{\P}{{\mathbb P}}
\newcommand{\R}{{\mathbb R}}
\newcommand{\T}{{\mathbb T}}

\def\br#1{\left(#1\right)}
\def\brb#1{\left[#1\right]}
\def\brs#1{\left\{#1\right\}}
\newcommand{\K}{K}

\theoremstyle{remark}
\newtheorem{rem}{Remark}

\theoremstyle{remark}
\newtheorem{notation}{Notation}

\theoremstyle{definition}
\newtheorem{definition}{Definition}

\theoremstyle{plain}
\newtheorem{theorem}{Theorem}[section]

\newtheorem{conjecture}{Conjecture}
\newtheorem{lemma}[theorem]{Lemma} 
\newtheorem{corollary}[theorem]{Corollary} 
\begin{document}

\title{Harmonic measure and SLE}
\date{}
\author{D.~Beliaev \and S.Smirnov}
\maketitle

\section{Introduction}

The motivation for this paper is twofold: to study multifractal spectrum of
the harmonic measure
and to better describe the geometry of Schramm's SLE curves. 
The main result of this paper is the Theorem \ref{th:spectrum} in which
we rigorously compute  the average spectrum of domains bounded by  SLE curves. 
Several results  can be easily
derived from this theorem: dimension estimates of the boundary of
SLE hulls, H\"older continuity of SLE Riemann maps, H\"older continuity
of SLE trace, and more. We also would
like to point out that SLE seems to be {\em the only model} where the
spectrum (even average) of harmonic measure is non-trivial and known explicitly.

\subsection{Integral means spectrum} \label{s:harmonic}

There are several equivalent definitions of harmonic measure that are
useful in different contexts. For a domain $\Omega$ with a regular boundary we
define the harmonic measure with a pole at $z\in \Omega$ as the exit
distribution of the standard Brownian motion started at $z$. Namely,
$\omega_{z}(A)=\P(B_\tau^z\in A)$, where $\tau=\inf\{t: B^{z}_t \not \in \Omega\}$
is the first time  the standard two-dimensional 
Brownian motion started at $z$ leaves $\Omega$.

Alternatively for a simply connected planar domain the harmonic
measure is the image of the normalized length on the unit circle
under the Riemann mapping that sends the origin to $z$.

It is easy to see that harmonic measure depends on $z$ in a smooth (actually
harmonic) way, thus the geometric properties do not depend on the choice of
the pole. So we fix the pole to be the origin or infinity and eliminate it
from notation.

Over the last twenty years it became clear  that many extremal problems in the geometric function
theory are related to the geometrical properties of harmonic measure
and the proper language for these problems is {\em the multifractal analysis}.

Multifractal analysis operates with different  spectra of measures
and relations between them. In this paper we study the harmonic measure on simply connected domains,
so we give the rigorous definition for this case only.

Let $\Omega=\C\setminus K$ where $K$ is a connected compact set and let
$\phi$ be a Riemann mapping from the complement of the unit disc $\D_-$ onto
$\Omega$ such that $\phi(\infty)=\infty$. {\em The integral means spectrum}
of $\phi$ (or $\Omega$) is defined as
$$
\beta_\phi(t)=\beta_\Omega(t)=\limsup_{r \to 1+}
\frac{\log \int |\phi(re^{i\theta})|^td\theta}{-\log(r-1)}.
$$
{\em The universal integral means spectrum} is defined as
$$
B(t)=\sup \beta_\Omega(t),
$$
where supremum is taken over all simply connected domains with compact boundary.

On the basis of work of Brennan,
Carleson, Clunie, Jones, Makarov, Pommerenke and
computer experiments for quadratic Julia sets
Kraetzer \cite{Kraetzer} in 1996 formulated the following universal conjecture:
\begin{eqnarray*}
B(t) = &t^2/4,& \quad |t|<2,\\
B(t) = &|t|-1,& \quad |t|\ge 2.
\end{eqnarray*}

It is known that many other conjectures follow from Kraetzer conjecture. In particular,
Brennan's conjecture \cite{Brennan} about integrability of $|\psi'|$ where $\psi$ is 
a conformal map from a domain to the unit disc is equivalent to $B(-2)=1$, while 
Carleson-Jones conjecture \cite{CaJo} about the decay rate of 
coefficient of a univalent
function and the growth rate of the length of the Green's lines is
 equivalent to $B(1)=1/4$.

There are many partial results in both directions: 
estimates of $B(t)$ from
above and below (see surveys \cite{BeSmECM,HeSo}).  
Upper bounds are more difficult and they still not that far from the trivial
bound $B(1)\le 1/2$. Currently the best upper bound is $B(1)\le 0.46$ \cite{HeSh}. 
Until recently lower bounds were also quite far from
the conjectured value. 

The main problem in finding lower bounds is that it is almost impossible
to compute the spectrum explicitly for any non-trivial domain. 
The origin of difficulties is easy  to see:  only fractal
domains have interesting spectrum, but for them the boundary
behavior of $|\phi'(re^{i\theta})|^t$ depends on $\theta$ in a very non
smooth way, making it hard to find the average growth rate.

We claim that in order to overcome these problems one should work with
regular random fractals instead of deterministic ones. For random fractals it
is natural to study {\em the average integral means spectrum} which is defined as
$$
\bar\beta(t)=
\limsup_{r\to 1}\frac{\log \int \E\brb{|\phi'(re^{i\theta})|^t}d\theta}{-\log|r-1|},
$$

The advantage of this approach it that for many random fractals the average boundary behavior
of $|\phi'|$ is a  very smooth function of $\theta$. 
Therefore it is sufficient to study average 
behavior along any particular radius. Regular
(random) fractals are invariant under some (random)
transformation, making $\E|\phi'|^t$   a solution of a specific equation.
Solving this equation one can find the
average spectrum.

Note that $\bar\beta(t)$ and $\beta(t)$ do not necessarily coincide. 
It can even happen (and in this paper we
consider exactly this case) that $\bar\beta(t)$ is  not a spectrum
of any particular domain. But Makarov's fractal approximation \cite{Makarov} implies that
$\bar\beta(t)\le B(t)$, so any function $\bar\beta$ gives a lower bound on $B$.

Another important notion is the {\em dimension} or 
{\em multifractal spectrum} of harmonic measure
 which can be non-rigorously defined as
$$
f(\alpha)=dim\{z: \omega( B(z,r)) \approx r^\alpha\},
$$
where $\omega (B(z,r))$ is the harmonic measure of the disc of radius $r$ centred at $z$.
There are several ways to make this definition rigorous, leading to slightly 
different notions of spectrum. But it is known \cite{Makarov} that the universal
spectrum $F(\alpha)=\sup f(\alpha)$ is the same for all definitions of $f(\alpha)$.

There is no simple relation between the integral means spectrum 
and the dimension spectrum for a given domain.
But for  regular (in some sense) fractals they are related by a Legendre type
transform. It is also known \cite{Makarov} that the universal spectra are related
by a Legendre type transform:
\begin{eqnarray*}
F(\alpha)&=&\inf_t(t+\alpha(B(t)+1-t)),\\
B(t)&=& \sup_{\alpha>0} \frac{F(\alpha)-t}{\alpha}+t-1.
\end{eqnarray*}

\subsection{Schramm-Loewner Evolution} \label{s:SLE}

It is a common belief that planar lattice models at criticality 
 have a conformally invariant scaling limits as
mesh of the lattice tends to zero. Schramm \cite{Schramm} 
introduced a one parametric family of random curves
which are called $SLE_\kappa$ (SLE stands for Stochastic 
Loewner Evolution or Schramm-Loewner Evolution) that are the only possible 
limits of cluster perimeters for critical lattice models.
It turned out to be also a very useful tool in many related problems.

In this section we give the definition of SLE and the necessary background.
The discussion of various versions of SLE and relations between them
can be found in Lawler's book \cite{Lawler05}. 

To define SLE we need a classical tool from complex analysis: the Loewner
evolution. In general this is a method to describe by an ODE the evolution of 
the Riemann map from a growing (shrinking) domain to a uniformization domain. 
In this paper we use the radial Loewner evolution (where uniformization 
domain is the complement of the unit disc) and its modifications.

\begin{definition}
The radial Loewner evolution in the complement of the unit disc with driving
function $\xi(t):\R_+\to\T$ is the solution of the following ODE
\begin{equation}
\partial_t g_t(z)=g_t(z)\frac{\xi(t)+g_t(z)}{\xi(t)-g_t(z)}, \qquad g_0(z)=z.
\label{eq:SLEdef}
\end{equation}
\end{definition}

It is a classical fact  \cite{Lawler05} that for any driving function $\xi$ $g_t$ is a conformal map
from $\Omega_t\to\D_-$ where $\D_-$ is the complement of 
the unit disc and $\Omega_t=\D_-\setminus\K_t$ is the set of all
points where solution of \eqref{eq:SLEdef} exists up to the time $t$.

The {\em Schramm-Loewner Evolution} $SLE_\kappa$ is defined as a Loewner
evolution driven by the Brownian motion with speed $\sqrt\kappa$ on the unit circle, namely
$\xi(t)=e^{i \sqrt{\kappa} B_t}$ where $B_t$ is the standard Brownian
motion and $\kappa$ is a positive parameter. 
Since $\xi$ is random, we obtain a family of random sets.
The corresponding family of compacts $K_t$ is also called SLE 
(or the {\em hull of SLE}).

A number of theorems was already established about SLE curves. Rohde and Schramm
\cite{RoSch} proved that SLE is a.s. generated by a curve. Namely,
almost surely there is a random curve 
$\gamma$ (called {\em trace}) such that $\Omega_t$ is
the unbounded component of $\D_- \setminus \gamma_t$,
where $\gamma_t=\gamma([0,t])$. The trace is almost 
surely a simple curve when $\kappa \le 4$. In this case 
the hull $K_t$ is the same as the curve $\gamma_t$. For $\kappa\ge 8$
the trace $\gamma_t$ is a space-filling curve. In the same paper they also
proved that almost surely the Minkowski (and hence  the Hausdorff)
dimension of the $SLE_\kappa$ trace is no more  than $1+\kappa/8$
for $\kappa\le 8$. Beffara \cite{Beffara1} proved that this estimate
is sharp for $\kappa=6$, later expanding the result to all $\kappa \le 8$.
Lind \cite{Lind} proved that the trace is
H\"older continuous. 

Another natural object is the {\em boundary} of  SLE hull, 
namely the boundary of $K_t$. For $\kappa\le 4$ the boundary 
of SLE is the same as SLE trace (since the trace is
a simple curve).  For $\kappa>4$ the boundary is the subset 
of the trace. Rohde and Schramm \cite{RoSch}
proved that for $\kappa>4$ the dimension of 
the boundary is no more than $1+2/\kappa$.

In \cite{Duplantier00, Duplantier03} physicist Duplantier 
using conformal field theory and quantum gravity methods predicted 
the average multifractal spectrum of SLE. He also conjectured  
that the boundary of $SLE_\kappa$ for $\kappa>4$ is in the same measure 
class as the trace of $SLE_{16/\kappa}$ 
(this is the so-called {\em duality} of SLE).

In this paper we rigorously compute the average integral means spectrum of SLE
and show that it coincides with  the Duplantier's prediction. 
This gives new proofs that dimension of the boundary is no 
more than $1+2/\kappa$ for $\kappa>4$ and
SLE maps are H\"older continuous,
and provides more evidence which supports the duality conjecture.

Since $\bar\beta$ is defined in the terms of a Riemann mapping, 
it is more convenient to
work with  $f_t=g_t^{-1}$. From the equation
(\ref{eq:SLEdef}) one can derive an equation on $f_t$. 
Unfortunately this equation involves
$f'_t$ as well as $\partial_t f_t$, so we have a PDE instead of ODE.

There is another approach which leads to a nice equation.
 Changing the direction of the
flow defined by the equation \eqref{eq:SLEdef} we get the equation for ``inverse'' function $g_{-t}$. 
For a given driving function $\xi$, maps $g_t^{-1}$ and $g_{-t}$ are
different, but in the case of Brownian motion they have the same 
distribution.The precise
meaning  is given by the following lemma
(which is an analog of the Lemma 3.1 from \cite{RoSch}):

\begin{lemma}
\label{lemma:timereverse}
Let $g_t$ be a radial SLE, then for all $t \in \R$ the map $z\mapsto g_{-t}(z)$
has the same distribution
as the map $z \mapsto \hat{f}_t(z)/\xi_t$, where $\hat{f}_t(z)=g_t^{-1}(z\xi_t)$.
\end{lemma}

\begin{proof}
Fix $s\in \R$. Let $\hat\xi(t)=\xi(s+t)/\xi(s)$. Then
$\hat\xi$ has  the same distribution as $\xi$.
Let
$$
\hat{g}_t(z)=g_{s+t}(g_{s}^{-1}(z \xi(s)))/\xi(s).
$$
It is easy to check that $\hat{g}_0(z)=z$ and
$$
\hat{g}_{-s}(z)=g_0( g_{s}^{-1}(z \xi(s)))/\xi(s)=
\hat{f}_{s}(z)/\xi(s).
$$
Differentiating $\hat{g}_{t}(z)$ with respect to $t$ we obtain
$$
\partial_t \hat{g}_t(z)
=\hat{g}_t(z)\frac{\hat\xi (t)+\hat{g}_t(z)}{\hat\xi (t)-\hat{g}_t(z)},
$$
hence $\hat g_t$ has the same distribution as SLE.
\end{proof}

This lemma proves that the solution of the equation
\begin{equation}
\partial_t f_t(z)=f_t(z)\frac{f_t(z)+\xi(t)}{f_t(z)-\xi(t)}, \qquad f_0(z)=z,
\label{eq:SLEdef2}
\end{equation}
where $\xi(t)=e^{i \sqrt{\kappa} B_t}$ has the same distribution as $g_t^{-1}$.
Abusing notations we call it also $SLE_\kappa$.

One of the most important properties of SLE is Markov property,
roughly speaking it means that the composition of two independent
copies of SLE is an SLE. The rigorous formulation is given by the following lemma.

\begin{lemma}
\label{lemma:Markov}
Let $f_\tau^{(1)}$ be an $SLE_\kappa$ driven by $\xi^{(1)}(\tau)$, $0<\tau<t$ and
$f_\tau^{(2)}$ be an $SLE_\kappa$ driven by $\xi^{(2)}(\tau)$, $0<\tau<s$,
where $\xi^{(1)}$ and $\xi^{(2)}$ are two independent Brownian motions on the circle.
Then $f_{s+t}(z)=f^{(2)}_t(f_t(z)/\xi^{(1)}(t))\xi^{(1)}(t)$ is $SLE_\kappa$ at time $t+s$.
\end{lemma}

\begin{proof}
This composition is the solution of Loewner Evolution driven by $\xi(\tau)$, where
\begin{equation*}
\xi(\tau)=\left\{
\begin{array}{ll}
\xi^{(1)}(\tau), & 0<\tau\le t, \\
\xi^{(2)}(\tau-t)\xi^{(1)}(t), & t<\tau\le t+s.
\end{array}
\right.
\end{equation*}
It is easy to see that $\xi(\tau)$ is also a Brownian motion on the circle with the same
speed $\sqrt{\kappa}$, hence $f_{t+s}$ is also $SLE_\kappa$.
\end{proof}

We will need yet another modification of SLE which is in fact a manifestation
of stationarity of radial SLE.

\begin{definition}
Let $\xi(t)=\exp(i\sqrt{\kappa}B_t)$ be a two-sided Brownian motion on the unit circle. The whole
plane $SLE_\kappa$ is the family of conformal maps $g_t$ satisfying
$$
\partial_t g_t(z)=g_t(z)\frac{\xi(t)+g_t(z)}{\xi(t)-g_t(z)},
$$
with initial condition
$$
\lim_{t\to-\infty}e^{t}g_t(z)=z, \qquad z\in \C\setminus \{0\}.
$$
\end{definition}

The whole-plane SLE satisfies the same differential equation as
the radial SLE, the difference is in the initial conditions. One
can think about the whole-plane SLE as about the radial SLE
started at $t=-\infty$. And this is the way to construct the
whole-plane SLE and prove the existence. Proposition 4.21 in
\cite{Lawler05} proves that the whole-plane Loewner Evolution $g_t$
with the driving function $\xi(t)$ is the limit as $s\to -\infty$ of
the following maps: $g^{(s)}_t(z)=e^{-t}z$ if $t\le s$,
$g^{(s)}_t(z)$ is the solution to (\ref{eq:SLEdef}) with initial
condition $g^{(s)}_s(x)=e^{-s}z$. The same is also true for inverse
maps.

We use this argument to prove that there is a limit of $e^{-t}f_t$ as $t\to \infty$.

\begin{lemma}
\label{lemma:stationarity}
Let $f_t$ be a radial $SLE_\kappa$ then there is a limit in the sense of distribution
of $e^{-t}f_t(z)$ as $t\to\infty$.
\end{lemma}

\begin{proof}
The function $e^{-t}f_t$ is exactly the function which is used to define the whole-plane
SLE. Multiplication by the exponent corresponds to the shift in time in the driving function.
The function $e^{-t}f_t(z)$ has the same distribution as the inverse of
$g^{(-t)}_0(z)$, hence it converges to $F_0$, where $F_\tau=g^{-1}_\tau$ and $g_\tau$
is a whole-plane SLE.
\end{proof}

\subsection{Results, conjectures, and organization of the paper}

It is easy to see that the geometry near ``the tip'' of SLE (the point of growth)
is different from the geometry near ``generic'' points. This means that 
for some problems it is more convenient to work with the so-called {\em bulk}
of SLE i.e. the part of the SLE hull which is away from the tip. 
In the following  theorem we compute the average spectrum of SLE hull and SLE bulk.   

\noindent{\bf Theorem \ref{th:spectru}}{\em \ \ 
The average integral means spectrum $\bar\beta(t)$ of SLE is equal to
$$
\begin{aligned}
-t+\kappa\frac{4+\kappa-\sqrt{(4+\kappa)^2-8t\kappa}}{4\kappa} & \qquad t \le -1-\frac{3\kappa}{8}, \\
-t+\frac{(4+\kappa)(4+\kappa-\sqrt{(4+\kappa)^2-8t\kappa}\ )}{4\kappa}  & \qquad  -1-\frac{3\kappa}{8}\le t \le \frac{3(4+\kappa)^2}{32\kappa}, \\
 t-\frac{(4+\kappa)^2}{16\kappa} & \qquad t\ge \frac{3(4+\kappa)^2}{32\kappa}.
\end{aligned}
$$
The average integral means spectrum $\bar\beta(t)$ of the  bulk of SLE is equal to
$$
\begin{aligned}
5
 -t+\frac{(4+\kappa)(4+\kappa-\sqrt{(4+\kappa)^2-8t\kappa}\ )}{4\kappa} & 
\qquad   t \le \frac{3(4+\kappa)^2}{32\kappa}, \\
 t-\frac{(4+\kappa)^2}{16\kappa}, & \qquad t\ge \frac{3(4+\kappa)^2}{32\kappa}.
\end{aligned}
$$
}

\begin{rem}
The local structure of the SLE bulk  is the same for all versions of SLE which means that 
they all have the same average spectrum.
\end{rem}

\begin{rem}
To prove this theorem we show that
$$
\E|f'(r e^{i\theta})|^t \asymp (r-1)^\beta ((r-1)^2+\theta^2)^\gamma,
$$ 
where $\beta$ and $\gamma$ are given by \eqref{eq:beta} and \eqref{eq:gamma}.
We would like to point out that $\beta$ and $\gamma$ are local exponents 
so they are the same for different versions of SLE.
\end{rem}

There are several corollaries that one can easily derive from Theorem \ref{th:spectru}

\begin{corollary}
SLE map $f$ is H\"older continuous with any exponent 
less than 
$$
\alpha_\kappa=1-\frac{1}{\mu}-\sqrt{\frac{1}{\mu^2}+\frac{2}{\mu}},
$$
where $\mu=(4+\kappa)^2/4\kappa$.
\end{corollary}

\begin{corollary}
The Hausdorff dimension of the boundary of the SLE hull for $\kappa\ge 4$ 
is at most $1+2/\kappa$.
\end{corollary}

\begin{corollary}
SLE trace with natural parametrization is H\"older continuous.
\end{corollary}

The first two results a conjectured to be sharp. They both have been 
previously published in \cite{Kang07} and \cite{RoSch} correspondingly. 
Both results can be easily derived from the properties of the spectrum 
(see \cite{Makarov})  and Theorem \ref{th:spectru}. 

The third corollary first appeared in a paper by Lind \cite{Lind}
where she uses derivatives estimates by Rohde and Schramm. One 
can use Theorem \ref{th:spectru} to prove this result.

The Theorem \ref{th:spectru} gives the average spectrum of SLE. The question 
about spectra of individual realizations of SLE remains open. 
We believe that with probability one they all have the same spectrum $\beta(t)$
which we call the a.s. spectrum. 

It is immediate that the tangent 
line at $t=3(4+\kappa)^2/32\kappa$ 
intersects $y$-axis at $-(4+\kappa)^2/16\kappa<-1$. 
This contradicts Makarov's characterization of possible 
spectra \cite{Makarov} which means that $\bar\beta$ 
can not be a spectrum of any given domain. In particular $\bar\beta$
is not the a.s. spectrum of SLE. On the other hand it suggests that 
the following conjecture is true.

\begin{conjecture}
\label{conj}
Let $t_{min}$ and $t_{max}$ be the two points such that the tangent to $\bar\beta(t)$ 
intersects the $y$-axis at $-1$. The almost sure value of the spectrum is equal to $\bar\beta(t)$
for $t_{min}\le t \le t_{max}$ and continues as the tangents for $t<t_{min}$ and $t>t_{max}$. Explicit 
formulas for $t_{min}$, $t_{max}$, and tangent lines are given in \eqref{eq:tminmax} and \eqref{eq:tangent}. See Figure \ref{pic:spectrum} for plots of $\beta$ and $\bar\beta$.
\end{conjecture}

The rest of the paper is organized in the following way. In the first part of the Section 
\ref{s:spectrum} we discuss Duplantier's prediction and the Conjecture \ref{conj}.
In the second part we compute the moments of 
$|f'|$ and prove Theorem \ref{th:spectru}.
In the Section \ref{s:levy} me make some remarks about
possible generalizations of SLE. In the last Section \ref{s:as} we explain 
a possible approach to the Conjecture \ref{conj}.

{\bf Acknowledgments. } Work supported in part by Swiss National Science foundation,
G\"oran Gustafsson Foundation, and Knut and Alice Wallenberg Foundation.

\section{Integral means spectrum of SLE}
\label{s:spectrum}

\subsection{Duplantier's prediction for the spectrum of the bulk}

In 2000 physicist Duplantier  predicted \cite{Duplantier00,Duplantier03} 
by the means of quantum gravity that
the Hausdorff dimension spectrum of the bulk of SLE  is
$$
f(\alpha)=\alpha-\frac{(25-c)(\alpha-1)^2}{12(2\alpha-1)},
$$
where $c$  is  the central charge which is related to $\kappa$ by
$$
c=\frac{(6-\kappa)(6-16/\kappa)}{4}.
$$
The negative values of $f$ do not have a simple geometric interpretation,
they correspond to negative dimensions (see papers by Mandelbrot 
\cite{Mandelbrot90,Mandelbrot03})
which appear only in the random setting. They correspond to the events
that have zero probability in the limit, but appear on finite scales as exceptional events.
There is another interpretation in terms of $\beta$ beta spectrum which we explain below.

Since negative values of $f$ correspond to zero robability events, 
it makes sense to introduce the positive part of the spectrum:  $f^+=\max\{f,0\}$.  
We believe  that $f^+$ is the almost sure value 
of the dimension spectrum. This is the dimension spectrum counterpart of Conjecture 
\ref{conj}. The function $f^+$ is equal to $f$ for  $\alpha\in [\alpha_{min},
\alpha_{max}]$, where
$$
\begin{aligned}
\alpha_{min}&=\frac{16+4\kappa+\kappa^2-2\sqrt{2}\sqrt{16\kappa+10\kappa^2+\kappa^3}}
{(4-\kappa)^2}, & \kappa\ne 4, \\
\alpha_{max}&=\frac{16+4\kappa+\kappa^2+2\sqrt{2}\sqrt{16\kappa+10\kappa^2+\kappa^3}}
{(4-\kappa)^2},  & \kappa\ne 4, \\
\alpha_{min}&=\frac{2}{3},  & \kappa=4,\\
\alpha_{max}&=\infty,  & \kappa=4.
\end{aligned}
$$
It is known (see \cite{Makarov}) that for regular fractals
$\beta(t)$ spectrum is related to $f(\alpha)$ spectrum by the Legendre transform. 
We believe those relations to hold for SLE as well:
$$
\begin{aligned}
\beta(t)-t+1=\sup_{\alpha>0}\,(f(\alpha)-t)/\alpha,\\
f(\alpha)=\inf_{t}\,(t+\alpha(\beta(t)-t+1)).
\end{aligned}
$$

The Legendre transform of $f^+$ is supposed to be equal to the almost sure value of the
integral means spectrum $\beta(t)$, while the Legendre transform of $f$ is  believed to be
equal to the average integral means spectrum $\bar\beta(t)$.

The Legendre transform of $f^+$ has two phase transitions: one for negative $t$ and one for
positive. The Legendre transform of $f^+$ is equal to
\begin{equation}
\begin{aligned}
\beta(t)&=t\br{1-\frac{1}{\alpha_{min}}}-1, &  t\le t_{min},
\\
\beta(t)&=-t+\frac{(4+\kappa)\br{4+\kappa-\sqrt{(4+\kappa)^2-8 t \kappa}}}{4\kappa},&
  t_{min}<t<t_{max},
\\
\beta(t)&= t\br{1-\frac{1}{\alpha_{max}}}-1, &  t\ge t_{max},
\end{aligned}
\label{eq:betaduplantier}
\end{equation}
where
$$
\begin{aligned}
t_{min}&=-f'(\alpha_{min})\alpha_{min},&\quad \kappa>0,\\
t_{max}&=-f'(\alpha_{max})\alpha_{max},& \quad \kappa\ne 4,\\
t_{max}&= 3/2, & \quad \kappa=4.
\end{aligned}
$$
We can also express $t_{min}$ and $t_{max}$ in terms of $\mu=4/\kappa+2+\kappa/4=(4+\kappa)^2/4\kappa$:
\begin{equation}
\begin{aligned}
t_{min}&=\frac{-1-2\mu-(1+\mu)\sqrt{1+2\mu}}{\mu},\\
t_{max}&=\frac{-1-2\mu+(1+\mu)\sqrt{1+2\mu}}{\mu}.
\end{aligned}
\label{eq:tminmax}
\end{equation}
And the linear functions in (\ref{eq:betaduplantier}) can be written
as
\begin{equation}
\begin{aligned}
\label{eq:tangent}
t \br{\frac{1}{\sqrt{1-{2 t_{min}}/{\mu}}}-1}-1 \\
t \br{\frac{1}{\sqrt{1-{2 t_{max}}/{\mu}}}-1}-1.
\end{aligned}
\end{equation}

For convenience we introduce
$$
\tilde\beta(t)=-t+\frac{(4+\kappa)\br{4+\kappa-\sqrt{(4+\kappa)^2-8 t \kappa}}}{4\kappa},
$$
which is the analytic part of the spectrum and defined for all $t<(4+\kappa)^2/8\kappa$.
This function is the analytic part of the Legendre transform of $f$. The
critical points $t_{max}$ and $t_{min}$ are the points where the tangent line
to the graph of $\bar\beta(t)$
intersects the $y$-axis at $-1$. The Legendre transform of $f^{+}$ is
equal to $\tilde \beta(t)$ between these two critical points and than continues as
a linear function.

Note that Makarov's theorem \cite{Makarov} states that all possible
integral means spectra  satisfy the following
conditions: they are non-negative convex functions bounded by the
universal spectrum such that tangent line at any point intersects $y$-axis
between $0$ and $-1$. So there is another way to describe the
Legendre transform of $f^{+}$: it coincides with $\tilde \beta$ as
long as this does not contradict Makarov's criteria and than
continues in the only possible way.

If we do not cut off the negative part of $f$, then the picture is a bit different.
There is no phase transition for negative $t$. For positive $t$ phase transition occurs
later, and it happens because the derivative of $f(\alpha)$ is bounded at infinity. For large $\alpha$
$$
f(\alpha)= \alpha\br{1-\frac{(4+\kappa)^2}{16\kappa}}+\frac{3(4+\kappa)^2}{32\kappa}+
O\br{\frac{1}{\alpha}},
$$
hence
$$
\begin{aligned}
\bar\beta(t)&=-t+\frac{(4+\kappa)\br{4+\kappa-\sqrt{(4+\kappa)^2-8 t \kappa}}}{4\kappa},&
\quad t\le\frac{3(4+\kappa)^2}{32\kappa},
\\
\bar\beta(t)&=1-\frac{(4+\kappa)^2}{16\kappa}+t-1=t-\frac{(4+\kappa)^2}{16\kappa}, &
\quad t>\frac{3(4+\kappa)^2}{32\kappa}.
\end{aligned}
$$

The explanation of this phase transition is rather simple. It is obvious that $\bar\beta(t)$ is
a convex function, and it follows from
Makarov's fractal approximation that the average spectrum is bounded by the universal spectrum.
It is known that for the large values of $|t|$ the universal spectrum is equal to $|t|-1$.
Altogether it implies that $|\bar\beta'(t)|\le 1$ and if it is equal to $1$ at some point
then $\bar\beta$ should be linear after this point. And $\bar\beta'=1$ exactly at
$t=3(4+\kappa)^2/32\kappa$.

\begin{figure}
    \centering
        \includegraphics[width=10cm]{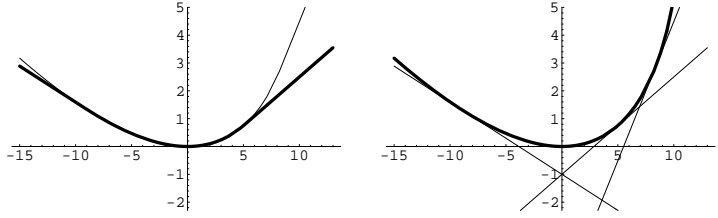}
    \caption{Plots of $\beta$ and $\bar\beta$ spectra.}
    \label{pic:spectrum}
\end{figure}

 The left part of the Figure \ref{pic:spectrum} shows plots of $\beta$ (thick line)
 and $\bar\beta$ for $\kappa=0.2$. On the right part one can see the plot of the
 $\tilde\beta$ (thick) and three tangent lines, two of them are crossing the
 $y$-axis at $-1$ and define $\beta$ after phase transition, the third one has
 a slope $1$ and defines $\bar\beta$ after the phase transition.

\subsection{Rigorous computation of the spectrum}

In this section we compute the average integral means spectrum of SLE (and its bulk) and
show that it coincides with the Legendre transform of the dimension spectrum
predicted by Duplantier.

The average integral means spectrum is the growth rate of
 $\tilde F(z,\tau)=\E\brb{|f'_\tau(z)|^t}$, where $f_\tau$ is
a radial $SLE_\kappa$. Actually, this function depends also on $t$ and $\kappa$, 
but they are fixed throughout the proof and we will not mention this dependence to simplify the notation.

\begin{lemma}
\label{lemma:martingale}
The function $\tilde F(z,\tau)$ is a solution of
\begin{equation}
\label{eq:martingale0}
\begin{aligned}
t \frac{r^4+4 r^2(1-r \cos \theta)-1}{(r^2-2 r \cos \theta +1)^2}\tilde F+
\frac{r(r^2-1)}{r^2- 2 r \cos\theta+1}\tilde F_r -
\\
\frac{2 r \sin\theta}{r^2-2 r \cos \theta +1} \tilde  F_\theta+
\frac{\kappa}{2} \tilde F_{\theta,\theta}-\tilde F_\tau=0.
\end{aligned}
\end{equation}
\end{lemma}

\begin{proof}
 The idea of the proof is to construct a martingale ${\mathcal M}_s$ 
 (w.r.t filtration defining SLE) which involves $\tilde F$.
 The $d s$ term in its  It\^o derivative should vanish. This will give us a partial
 differential equation on $\tilde F$.
We set
$$
{\mathcal M}_s=\E \brb{ |f_{\tau}'(z)|^t \mid {\mathcal F}_s }.
$$
By  the Lemma \ref{lemma:Markov}
\begin{eqnarray*}
\E \brb{ |f_{\tau}'(z)|^t \mid {\mathcal F}_s }&=&
\E \brb{ | f_s'(z)|^t |f_{\tau-s}'( f_s(z)/\xi_s)|^t \mid {\mathcal F}_s }
\\
&=&| f_s'(z)|^t\tilde F(z_s,\tau-s),
\end{eqnarray*}
where $z_s=f_s(z)/\xi_s$.

We will need derivatives of $z_s$ and $| f_s'|^t$
\begin{eqnarray*}
\partial_t \log| f_s'(z)|&=&
\Re \frac{\partial_z f_s \frac{ f_s+\xi_s}{ f_s-\xi_s}}{ f_s'}=
\Re\brb{\frac{ f_s +\xi_s}{ f_s -\xi_s} - \frac{2 \xi_s  f_s}{( f_s-\xi_s)^2}}
\\ & =&
\Re \frac{z_s^2-1-2 z_s}{(z_s-1)^2} =
\frac{r^4+4 r^2(1-r \cos \theta)-1}{(r^2-2 r \cos \theta +1)^2},
\end{eqnarray*}
where $z_s=r\exp(i \theta)$.
Next we have to find the derivative of $z_s=r e^{i\theta}$
\begin{equation*}
 d \log z_s =d \log r + i d \theta= d \log  f_s -  i\sqrt{\kappa} d B_s,
\end{equation*}
where
\begin{equation*}
d \log  f_s =\frac{d  f_s}{f_s}=\frac{z_s+1}{z_s-1}d s.
\end{equation*}
Writing everything in terms of $r$ and $\theta$ we get
$$
\begin{aligned}
d \log r +i d \theta = \frac{z_s+1}{z_s-1}d s - i\sqrt{\kappa} d B_s =
\\
\frac{r^2-1}{r^2- 2 r \cos\theta+1}d s +
i \br{ -\frac{2 r \sin\theta}{r^2-2 r \cos \theta +1}d s - \sqrt{\kappa}d B_s}.
\end{aligned}
$$
Summing it all up we obtain
\begin{eqnarray}
\partial_t \log| f_s'(z)| = \frac{r^4+4 r^2(1-r \cos \theta)-1}{(r^2-2 r \cos \theta +1)^2},
\\
d \theta = -\frac{2 r \sin\theta}{r^2-2 r \cos \theta +1}d s - \sqrt{\kappa}d B_s,
\\
d r = r d \log r =\frac{r(r^2-1)}{r^2- 2 r \cos\theta+1}d s .
\end{eqnarray}
Let us write $F(z,\tau)$ as  $F(r,\theta, \tau)$.
The $d s$ term in It\^o derivative of ${\mathcal M}$ is equal to
\begin{eqnarray*}
| f_s'(z)|^t
\left(
t\frac{r^4+4 r^2(1-r \cos \theta)-1}{(r^2-2 r \cos \theta +1)^2}\tilde F+
\frac{r(r^2-1)}{r^2- 2 r \cos\theta+1}\tilde F_r - \right.
\\
\left. \frac{2 r \sin\theta}{r^2-2 r \cos \theta +1}  \tilde F_\theta+
\frac{\kappa}{2} \tilde F_{\theta,\theta}-\tilde F_\tau\right).
\end{eqnarray*}
This derivative should be $0$ and, since $f_s$ is a univalent  function and its derivative
never vanishes,  $\tilde F$ is a solution of (\ref{eq:martingale0}).
\end{proof}

By the Lemma \ref{lemma:stationarity} there is a limit of $e^{-\tau}f_\tau$ as $\tau\to\infty$.
Hence we can introduce
$$
F(z)=\lim_{\tau\to\infty} e^{-\tau t}F(z,\tau).
$$
Passing to the limit in (\ref{eq:martingale0}) we can see that $F(z)$ is a solution of

\begin{equation}
\label{eq:martingale}
\begin{aligned}
t \br{\frac{r^4+4 r^2(1-r \cos \theta)-1}{(r^2-2 r \cos \theta +1)^2}-1} F+
\frac{r(r^2-1)}{r^2- 2 r \cos\theta+1} F_r
\\
-\frac{2 r \sin\theta}{r^2-2 r \cos \theta +1}   F_\theta+
\frac{\kappa}{2}  F_{\theta,\theta}=0.
\end{aligned}
\end{equation}

\begin{notation}
We define two constants $\beta$ and $\gamma$:
\begin{eqnarray}
\gamma=\gamma(t,\kappa)=\frac{4+\kappa - \sqrt{(4+\kappa)^2-8 t \kappa}}{2\kappa},
\label{eq:gamma}
\\
\beta=\beta(t,\kappa)=t-\frac{(4+\kappa)\gamma}{2}.
\label{eq:beta}
\end{eqnarray}
It is easy to see that the second constant $\beta$ is equal to $-\tilde\beta$.
\end{notation}

Let us explain where these constants  come from. Roughly speaking
spectrum $\beta(t)$ is the growth rate of $F$ as $r\to 1$. $F$ is a
solution of the parabolic equation (\ref{eq:martingale}) that has a
singularities  at $|z|=1$ and $z=1$. Let us assume that $F$ has a
power series expansion near $1$. Then we can write power series
expansion of coefficients of (\ref{eq:martingale}) and assuming that
the leading term is $(r-1)^\beta((r-1)^2+\theta^2)^\gamma$ we get an
equation on $\beta$ and $\gamma$. Constants $\gamma$ and $\beta$ are
solution of these equation. Now let us explain why it makes sense to
consider this expansion. Instead of radial/whole-plane SLE we can
write the same martingale for the chordal SLE. In this case (if we
forget dependence on $\tau$) the equation will be
\begin{equation}
2t\frac{x^2-y^2}{(x^2+y^2)^2}F-\frac{2x}{x^2+y^2}F_x+\frac{2y}{x^2+y^2}F_y+\frac{\kappa}{2}F_{xx}=0.
\label{eq:chordal}
\end{equation}
This equation is ``tangent'' to \eqref{eq:martingale} at $r=1$ and $\theta=0$.

This equation has a solution of the form $y^\beta(x^2+y^2)^\gamma$,
where $\beta$ and $\gamma$ as above. Actually this is the way how we
found this exponents. This approach seems to be easier, but there
are two major problems. First it is not easy to argue that we can
neglect the derivative with respect to $\tau$. Another problem is
that $y^\beta(x^2+y^2)^\gamma$ can not be equal to $F$ since it
blows up at infinity and we have to show that the local behavior
does not depend on the boundary conditions at infinity.

Similar equation appeared in \cite{RoSch} and when this work was finished we learned from I. Gruzberg that
equation (\ref{eq:chordal}) appeared several years ago in the paper
by Hastings \cite{Hastings}.

\begin{theorem}
\label{th:spectrum}
Let
$$
t\le \frac{3(4+\kappa)^2}{32\kappa}.
$$
Then we have
$$
\E\brb{\int_{|z|=r}|F'_0(r e^{i\theta})|^td \theta}\asymp \br{\frac{1}{r-1}}^{\bar\beta(t)} ,
$$
where the expectation is taken for a;; whole-plane SLE maps $F_0=\lim e^{-\tau}f_\tau$ 
and $\bar\beta(t)$ is equal to
\begin{equation}
\label{eq:half}
\begin{aligned}
-\beta(t,\kappa),  &
\qquad t>-1-\frac{3\kappa}{8},\\
 -\beta(t,\kappa)-2\gamma(t,\kappa)-1,   &
 \qquad t\le -1-\frac{3\kappa}{8}.
\end{aligned}
\end{equation}

\end{theorem}

\begin{proof}
Suppose that we can find  positive functions  bounded away from the unit circle 
$\phi_+$ and $\phi_-$ such that
$\Lambda \phi_-<0$ and $\Lambda \phi_+>0$ then by the maximum principle any 
positive solution of \eqref{eq:martingale} is between $c_+\phi_+$ and $c_-\phi_-$,
where $c_+$ and $c_-$ are positive constants.

In the Lemma \ref{boundarysolution} we will construct such
functions $\phi_-$ and $\phi_+$.   They are of the form
$$
\phi_\pm =(r-1)^\beta(r^2-2r\cos\theta+1)^\gamma(-\log(r-1))^{\mp 1} g(r^2-2r\cos\theta+1),
$$
where $g>0$ for $r=1$. Both functions have the same polynomial growth rate as $r\to 1$, thus
$F$ has also the same growth rate. By Tonelli theorem
$$
\E\brb{\int|f'_\tau|^t}=\int \E\brb{|f'_\tau(r,\theta)|^t}d\theta \approx
\int(r-1)^\beta(r^2-2r\cos\theta+1)^\gamma d\theta,
$$
where $\approx$ means that functions have the same polynomial growth rate.
For $\gamma>-1/2$ the weight $(r^2-2r\cos\theta+1)^\gamma$ is integrable
up to the boundary and we immediately get
$$
\E\brb{\int_{|z|=r}|f'_\tau|^t}\approx \br{\frac{1}{r-1}}^{-\beta}.
$$
For $\gamma\le-1/2$ the situation is a bit different. In this case the integral of
the weight blows up as $(r-1)^{2\gamma+1}$. Which gives us
$\E\brb{\int|f_\tau'|^td\theta}\approx(r-1)^{\beta+2\gamma+1}$.
It is easy to check that $\gamma\le -1/2$ if and  only if $t\le -1-3\kappa/8$.
\end{proof}

\begin{rem}
The growth rate of $\E\brb{\int|f'|^t}$ is similar to $\bar\beta(t)$ predicted by Duplantier.
The phase transition at $t=-1-3\kappa/8$ is due to the exceptional behavior of SLE at the tip. If we
integrate over values of $\theta$ bounded away from $0$ then the weight $|z-1|^{2\gamma}$ does not blow up
and we have no phase transition at $t=-1-3\kappa/8$ any more. This gives us the spectrum of the bulk of SLE.
\end{rem}

We can also state this theorem in terms of average integral means
spectrum defined in the introduction. This
theorem proves that Duplantier's prediction for $\bar \beta(t)$ is correct.

\begin{theorem}
\label{th:spectru}

The average integral means spectrum $\bar\beta(t)$ of SLE is equal to
$$
\begin{aligned}
-t+\kappa\frac{4+\kappa-\sqrt{(4+\kappa)^2-8t\kappa}}{4\kappa} & \qquad t \le -1-\frac{3\kappa}{8}, \\
-t+\frac{(4+\kappa)(4+\kappa-\sqrt{(4+\kappa)^2-8t\kappa}\ )}{4\kappa}  & \qquad  -1-\frac{3\kappa}{8}\le t \le \frac{3(4+\kappa)^2}{32\kappa}, \\
 t-\frac{(4+\kappa)^2}{16\kappa} & \qquad t\ge \frac{3(4+\kappa)^2}{32\kappa}.
\end{aligned}
$$
The average integral means spectrum $\bar\beta(t)$ of the  bulk of SLE is equal to
$$
\begin{aligned}
 -t+\frac{(4+\kappa)(4+\kappa-\sqrt{(4+\kappa)^2-8t\kappa}\ )}{4\kappa} & 
\qquad   t \le \frac{3(4+\kappa)^2}{32\kappa}, \\
 t-\frac{(4+\kappa)^2}{16\kappa}, & \qquad t\ge \frac{3(4+\kappa)^2}{32\kappa}.
\end{aligned}
$$

\end{theorem}
\begin{proof}
The Theorem \ref{th:spectrum} gives us the value of $\bar\beta(t)$ for $t\le 3(4+\kappa)^2/32\kappa$.
Direct computations show that derivative of $-\beta(t,\kappa)$ at $t=3(4+\kappa)^2/32\kappa$
is equal to one. As we mentioned before, the $\bar\beta$ spectrum is a convex function bounded by
the universal spectrum, and the universal spectrum is equal to $|t|-1$ for the large values of $|t|$.
This means that if $\bar\beta'=1$ at some point then it should continue as a linear function with
slope one. Hence  $\bar\beta$ should continue as   $t-{(4+\kappa)^2}/{16\kappa}$
for $t>3(4+\kappa)^2/32\kappa$. Plugging in the values of $\beta$ and 
$\gamma$ we finish the proof of the theorem.
\end{proof}

To complete the proof of the Theorem \ref{th:spectrum} we have to construct functions $\phi_-$ and
$\phi_+$. We do it in three  steps, first we write the restriction of the equation
(\ref{eq:martingale}) to the unit circle, then we find a positive solution $g$ of the resulting equation.
Finally we construct $\phi_-$ and $\phi_+$ out of $g$.

We look for a solution in the following form:
$$
f(r,\theta)=(r-1)^\beta(r^2-2 r \cos\theta+1)^\gamma g(r^2-2 r \cos\theta +1).
$$
Plugging $f$ into (\ref{eq:martingale}), factoring $(r-1)^\beta (r^2-2 r \cos\theta+1)^\gamma$ out,
and taking $r=1$ we obtain a differential equation on $g(2-2\cos\theta)$
\begin{equation}
\begin{aligned}
(-2 t +4 \beta-2\gamma-2\gamma\kappa+\gamma^2\kappa+
(4 t - 4\beta +2\kappa \gamma)\cos\theta
\\-(2 t +\gamma(\gamma\kappa-2)\cos (2\theta))g(2-2 \cos\theta)
\\
+(2-2\cos\theta)(-2-\kappa+2\gamma\kappa+2\kappa\cos\theta-(\kappa+2\gamma\kappa-2)\cos(2\theta))
\\
\times g'(2-2\cos\theta)+ 2\kappa(2-2\cos\theta)(\sin\theta)^2 g''(2-2\cos\theta)=0.
\label{eq:boundary0}
\end{aligned}
\end{equation}

\begin{lemma}
The equation (\ref{eq:boundary0}) has a smooth (with possible exception at $\theta=0$)
positive bounded solution on the circle if and only if
\begin{equation}
\label{cond}
t\le \frac{3(4+\kappa)^2}{32\kappa}.
\end{equation}
\label{lemma:boundary}
\end{lemma}

\begin{proof}
Changing the variable to $x=2-2\cos\theta$ we rewrite (\ref{eq:boundary0}) as a
hypergeometric equation
\begin{equation}
\label{eq:boundary}
\begin{aligned}
\gamma(2+\kappa)g(x)+(8-2 x +\kappa(x-2)+2\gamma\kappa(x-4))g'(x)+ \\
\kappa(x-4)x g''(x)=0,
\end{aligned}
\end{equation}
which  has two independent solutions
$$
g_1(x)={_2F_1}(a,b,\frac{1}{2}+a+b,\frac{x}{4})
$$
and
$$
g_2(x)=x^{1/2-a-b}{_2F_1}(\frac{1}{2}-a,\frac{1}{2}-b,\frac{3}{2}-a-b,\frac{x}{4}),
$$
where
$$
\begin{aligned}
a=\gamma-\frac{1}{\kappa}-\frac{\sqrt{1-2t\kappa}}{\kappa},\\
b=\gamma-\frac{1}{\kappa}+\frac{\sqrt{1-2t\kappa}}{\kappa}.
\end{aligned}
$$
Function $g(2-2\cos\theta)$ should have a second derivative everywhere
on the unit circle except at the point $\theta=0$. This
means that $g(x)$ should have expansion $c+O(4-x)$ at the endpoint $4$.

Any solution of (\ref{eq:boundary0}) is a linear combination of 
$g_1$ and $g_2$:
$g=c_1 g_1+c_2 g_2$. We want to find coefficients $c_1$ and $c_2$ such that
this sum is bounded and has a correct expansion at $x=4$.

Expansions of $g_1$ and $g_2$ at $4$ are
\begin{eqnarray*}
g_1(x)&=&\frac{\sqrt{\pi}\Gamma(1/2+a+b)}{\Gamma(1/2+a)\Gamma(1/2+b)}-
\\
&&
\frac{\sqrt{\pi}\Gamma(1/2+a+b)}{\Gamma(a)\Gamma(b)}\sqrt{4-x}+
O(4-x),
\end{eqnarray*}
and
\begin{eqnarray*}
g_2(x)&=&\frac{2^{1-2a-2b}\sqrt{\pi}\Gamma(3/2-a-b)}{\Gamma(1-a)\Gamma(1-b)}-
\\
&&
\frac{2^{1-2a-2b}\sqrt{\pi}\Gamma(3/2-a-b)}{\Gamma(1/2-a)\Gamma(1/2-b)}\sqrt{4-x}+O(4-x).
\end{eqnarray*}

If $c_2\ne 0$ then $1/2-a-b$ should be nonnegative, otherwise $g$ is not bounded at $0$.
Note that
$$
\frac{1}{2} -a-b=\frac{4+\kappa-4\gamma\kappa}{2\kappa}
$$
which is nonnegative if and only if
$$
t\le \frac{3(4+\kappa)^2}{32\kappa}
$$
which is exactly the restriction from the statement of the lemma.
If $t> 3(4+\kappa)^2/32\kappa$ then $c_2=0$.  In this case  $g$ has
a correct expansion at $4$ if and only if $\Gamma(a)=0$ or $\Gamma(b)=0$, but
$1-2t\kappa<0$ so both $a$ and $b$ are not real number and gamma function has only real
roots.

We can introduce
$$
C=\frac{\Gamma(1/2+a+b)\Gamma(1/2-a)\Gamma(1/2-b)}{2^{1-2a-2b}\Gamma(a)\Gamma(b)\Gamma(3/2-a-b)},
$$
and
$$
g_3(x)=g_1(x)-C g_2(x).
$$
By construction $g_3(x)=\mathrm{const}+O(4-x)$ near $4$. Finally we have to prove that
$g_3$ is a positive function. Note that in (\ref{eq:boundary})
 $g$ and $g''$ have coefficients of different signs. Obviously, $g_3(0)=1$. Suppose that
 $g_3$ has a local minimum inside the interval $(0,4)$, then $g_3'=0$ and $g_3''\ge 0$ at this point,
 hence $g_3$ is also positive. Thus it is sufficient to check that $g_3(4)>0$. The value
of $g_3(4)$ is easy to evaluate
\begin{eqnarray*}
g_3(4)&=& \sqrt{\pi}\Gamma(1/2+a+b)
\\
&&\times
\br{\frac{1}{\Gamma(1/2+a)\Gamma(1/2+b)}-
\frac{\Gamma(1/2-a)\Gamma(1/2-b)}{\Gamma(a)\Gamma(b)\Gamma(1-a)\Gamma(1-b)}}
\\
&=&
\frac{\sqrt{\pi}\Gamma(1/2+a+b)\cos(\pi(a+b))}{\Gamma(1/2+a)\Gamma(1/2+b)\cos(\pi a)\cos(\pi b)}
\\
&=&
\pi^{-3/2}\Gamma(1/2+a+b)\cos(\pi(a+b))\Gamma(1/2-a)\Gamma(1/2-b).
\end{eqnarray*}
By (\ref{cond}), $a+b<1/2$, hence $\Gamma(1/2+a+b)\cos(\pi(a+b))>0$. Finally
we have to show that $\Gamma(1/2-a)\Gamma(1/2-b)>0$. We  consider two
different cases: when $t\le 1/2\kappa$ and $t >1/2\kappa$. In the second case
$a$ and $b$ are conjugated and $\Gamma(1/2-a)\Gamma(1/2-b)=|\Gamma(1/2-a)|^2>0$.
In the first case we will prove that $1/2-a>0$ and $1/2-b>0$. It is easy to see that
$1/2-b<1/2-a$, hence it is sufficient to prove that $1/2-b>0$. Recall that
$$
\frac{1}{2}-b=\frac{1}{2}-\gamma+\frac{1}{\kappa}-\frac{\sqrt{1-2t\kappa}}{\kappa},
$$
hence
$$
\partial_t(1/2-b)=\frac{1}{\sqrt{1-2\kappa t}}-\frac{2}{\sqrt{(4+\kappa)^2-8t\kappa}}>0.
$$
This means that $1/2-b$ has a minimum when $t=0$, this minimum is
$$
\frac{1}{2}-b(0)=\frac{1}{2}-\gamma(0)=\frac{1}{2}>0.
$$
This proves that $g_3(x)>0$ on $[0,4]$.
\end{proof}

\begin{lemma}
\label{boundarysolution}
Let $g$ be a positive  bounded solution of (\ref{eq:boundary0}) and
$$
\begin{aligned}
F=&f(r,\theta)(-\log(r-1))^\delta\\=
&(r-1)^\beta(r^2-2 r \cos\theta +1)^\gamma g(r^2-2 r \cos\theta +1)(-\log(r-1))^\delta.&
\end{aligned}
$$
Then
\begin{eqnarray*}
\Lambda F >0, \quad \delta<0,
\\
\Lambda F<0, \quad \delta>0,
\end{eqnarray*}
for $r$ sufficiently close to $1$.
\end{lemma}

\begin{proof}
Applying $\Lambda$ we find
$$
\Lambda F=(-\log(r-1))^\delta
\br{ \Lambda f -f \frac{r(r+1)\delta}{(r^2-2 r \cos \theta +1)(-\log(r-1))}}.
$$
By Lemma \ref{lemma:boundary} $\Lambda f=(r-1)^\beta(r^2-2r\cos\theta +1)^\gamma O(r-1)$,
hence
$$
\begin{aligned}
\Lambda F= &(-\log(r-1))^\delta(r-1)^\beta(r^2-2r\cos\theta +1)^\gamma
\\
&\times\br{O(r-1)- \frac{r(r+1)\delta(g(2-2\cos\theta)+O(r-1))}{w(-\log(r-1))}}.
\end{aligned}
$$
The sign of the main term is opposite to the sign of  $\delta$. This proves the claim.
\end{proof}

\begin{rem}
Note that we proved a  stronger result than announced in Theorem \ref{th:spectrum}:
$\E \int |F'|^t $ has growth rate $(r-1)^\beta$ up to a factor $\log^\delta(r-1)$
for {\em arbitrary small} $|\delta|$.
\end{rem}

\section{Loewner Evolution driven by other processes}
\label{s:levy}

It is known that Loewner Evolution can be defined for a very large
class of driving functions. In particular, they do not have to be
continuous. In  \cite{BeSmECM} we proposed to study L\'evy-Loewner
Evolution ($LLE$), which is the Loewner Evolution driven by a L\'evy
process (i.e. process with independent stationary increments). This
defines a very rich class of random fractals. It seems that it is
still possible to find the spectrum of harmonic measure for this
class explicitly.

In the fundamental Lemma \ref{lemma:martingale} we  only use  the
fact that the  Brownian motion is a L\'evy process. So the same
argument can be applied for $LLE$. As the result we get that
$F=\E\brb{|e^{-\tau}f_\tau'(z)|^t}$ is the solution of
$$
\begin{aligned}
t \br{\frac{r^4+4 r^2(1-r \cos \theta)-1}{(r^2-2 r \cos \theta +1)^2}-1} F+
\frac{r(r^2-1)}{r^2- 2 r \cos\theta+1} F_r
\\
-\frac{2 r \sin\theta}{r^2-2 r \cos \theta +1}   F_\theta+\Lambda F=0.
\end{aligned}
$$
where $\Lambda$ is the generator of the driving L\'evy process.
Thus again finding spectrum boils down to the analysis of a
parabolic type integro-differential equation. We have a freedom to
chose the driving process (and the generator $\Lambda$), so it seems
possible to find such driving process that this equation could be
solved and gives  large spectrum.

This paper was in preparation for a long time. During this period 
there appeared several paper studying the most natural $LLE$ where 
the driving force is a symmetric $\alpha$ stable process (or a sum of a 
Brownian motion and stable process). Unpublished computer experiments 
by Meyer \cite{Meyer} suggested that the spectrum for $1$-stable 
process could be large (and possibly  equal to the conjectured
universal spectrum). Unfortunately later work by
Gruzberg, Guan, Kadanoff, Oikonomou, Rohde, Rushkin, Winkel, 
and others \cite{ROKG,GuWi,Guan}
showed that this is wrong. But there is still a possibility that 
computer experiments exposed an existing phenomenon. It could be that
the integral means grow fast for a few (relatively) large scales and
when we approach the boundary their growth slows down. If this is true, 
one can use $LLE$ as a building block in a snowflake (or any other construction
which allows to replicate scales). 
In this way one can hope to construct a domain with large 
integral means on {\em all} scales.

\section{Almost sure value of the spectrum}
\label{s:as}

In this section we speculate about what should be done to 
prove that the almost sure value of the spectrum is given by \eqref{eq:betaduplantier}.

Let us introduce random variables $X_k{(n)}=|f'((1+2^{-n})e^{2\pi i k/2^n})|^t$.
The spectrum is the growth rate of $2^{-n}\sum_k X_k$. We know that 
$$
2^{-n}\sum_{k=1}^{2^n} \E X_k \asymp 2^{n\beta(t)}.
$$ 
We want to show that the probability 
\begin{equation}
\label{eq:probability}
\P\brs{2^{-n}|\sum X_k-\E X_k|>2^{n (\beta(t)-\delta)}}
\end{equation}
is summable for some positive $\delta$. This will clearly imply that 
spectrum of SLE is equal to $\beta(t)$ with probability one.

Conformal field theory considerations suggest 
that $X_k$ and $X_l$
are essentially independent if $|k-l|\gg 1$ 
(in other words the distance between points should be much larger 
than their distance to the boundary).   In fact it is believed that 
derivatives are essentially  independent if the distance between points is greater 
that any power (less than one) of the distance to the boundary. Let 
us exaggerate it a little bit more and assume that $X_k$ and $X_l$ 
are independent for any $k \ne l$.

Let us denote $X_k-\E X_k$ by $Y_k$. By Chebyshev inequality the probability \eqref{eq:probability} is less than
$$
\frac{\E|\sum Y_k|^{1+\epsilon}}{2^{n(1+\epsilon)(\beta(t)+1-\delta)}}.
$$
It is known (see \cite{BE}) that for independent random variables with zero mean 
$\E|\sum Y_k|^{1+\epsilon}\le c \sum \E |Y_k|^{1+\epsilon}$, where $c$ is an absolute constant
which does not depend on the number of terms. Using this 
we can estimate the fraction above by
\begin{equation}
\begin{aligned}
\frac{\sum \E |Y_k|^{1+\epsilon}}{2^{n(1+\epsilon)(\beta(t)+1-\delta)}}\le 
c\frac{2^n 2^{n\beta(t+t\epsilon)}}{2^{n(1+\epsilon)(\beta(t)+1-\delta)}}= 
\\
c 2^{n(1+\beta(t+t\epsilon)-\beta(t)-1+\delta-\epsilon\beta(t)-\epsilon+\epsilon\delta)}.
\end{aligned}
\label{eq:exp}
\end{equation}

For small $\epsilon<\epsilon_0(t)$ the exponent in the last formula is bounded by
$$
\begin{aligned}
n(\beta'(t)t\epsilon+\epsilon^{3/2}+\delta-\epsilon\beta(t)-\epsilon+\epsilon\delta)= \\
n(\epsilon(\beta'(t)t-\beta(t)-1)+\epsilon^{3/2}+\delta+\epsilon\delta).
\end{aligned}
$$
If $\beta'(t)t-\beta(t)-1=c(t)<0$, then we can find a small $\epsilon_t$ (depending on $t$ only)
 such that 
$\epsilon_t(\beta'(t)t-\beta(t)-1)+\epsilon_t^{3/2}<c(t)\epsilon_t/2$. Fix $\delta=-\epsilon_t c(t)/4$, then the exponent in \eqref{eq:exp} is negative. This implies that the probability in \eqref{eq:probability} is summable if $-1<\beta(t)-t\beta'(t)$. The last inequality means that the tangent line to $\beta$ at point $t$ intersects the $y$ axis above $-1$. This is exactly the condition which appeared in \eqref{eq:betaduplantier}. 

Thus, assuming the independence of derivatives, we can prove that the almost sure value of the spectrum is equal to $\bar\beta(t)$ for $t_{min}<t<t_{max}$. For other values of $t$ Makarov's theorem implies that the spectrum should continue as a straight line tangent to $\bar\beta(t)$ at $t_{min}$ and $t_{max}$ correspondingly.

\bibliography{sle}
\bibliographystyle{abbrv}

\end{document}